\documentclass[a4paper,reqno,12pt]{amsart}
\usepackage[margin=2.5cm,papersize={19.75cm,28cm}]{geometry}

\usepackage{hyperref}
\usepackage{amssymb}
\usepackage{eucal}
\usepackage{eufrak}
\usepackage{graphics}
\usepackage{graphicx}
\usepackage{epsfig}
\usepackage{multicol}
\usepackage{xypic}
\usepackage{amsmath}
\usepackage{wasysym}
\usepackage{dsfont}
\usepackage{xypic}
\usepackage{amsmath}
\usepackage{color}
\usepackage{tikz}
\usepackage{xy}
\xyoption{all}

\usetikzlibrary{matrix}

\newtheorem{theorem}{Theorem}[section]
\newtheorem{proposition}{Proposition}[section]
\newtheorem{definition}{Definition}[section]

\newtheorem{example}[definition]{Example}

\newtheorem{remarkth}[definition]{Remark}
\newenvironment{remark}{\begin{remarkth}\upshape}{\end{remarkth}}

\newcommand{\Flder}{\rightarrow}

\newcommand{\lc}{\left\{}
\newcommand{\rc}{\right\}}

\newcommand{\ra}{\rightarrow}

\newcommand{\bra}{\langle}
\newcommand{\ket}{\rangle}
\newcommand{\R}{\mathds{R}}      

\newcommand{\I}{\mathds{I}}
\newcommand{\lcf}{\lbrack\! \lbrack}
\newcommand{\rcf}{\rbrack\! \rbrack}

\usepackage{amssymb}

\newcommand{\proa}{A^*G \mbox{$\;$}_{\tau^*} \kern-3pt\times_\alpha
G \mbox{$\;$}_\beta \kern-3pt\times_{\tau^*} A^*G}

\newcommand{\ab}{\lbrack\!\lbrack \;, \; \rbrack\!\rbrack}
\def\lcf{\lbrack\! \lbrack}
\def\rcf{\rbrack\! \rbrack}

\begin{document}

\title[Optimal control of nonholonomic mechanical systems]{ A geometric 
approach to the optimal control of nonholonomic mechanical systems}

\author[A. Bloch]{Anthony Bloch}
\address{A. Bloch: Department of Mathematics, University of Michigan, Ann Arbor, MI 48109, USA} \email{abloch@umich.com}

\author[L. Colombo]{Leonardo Colombo}
\address{L. Colombo: Department of Mathematics, University of Michigan, Ann Arbor, MI 48109, USA} \email{ljcolomb@umich.edu}

\author[R. Gupta]{Rohit Gupta}
\address{R. Gupta: Department of Aerospace Engineering, University of Michigan, Ann Arbor, MI 48109,
USA}\email{rohitgpt@umich.edu}

\author[D. Mart\'in de Diego]{David Mart\'in de Diego}
\address{D.\ Mart\'{\i}n de Diego: Instituto de Ciencias Matem\'aticas (CSIC-UAM-UC3M-UCM), Calle Nicol\'as Cabrera 15, 28049 Madrid, Spain} \email{david.martin@icmat.es}



\thanks{\noindent This work has been partially supported by Grants   MTM2013-42870-P, MTM2009-08166-E, IRSES-project ``Geomech-246981'' and NSF grant INSPIRE-1363720 and NSF grant 1207693}

\begin{abstract}
In this paper, we describe a constrained Lagrangian and Hamiltonian
formalism for the optimal control of nonholonomic mechanical
systems. In particular, we aim to minimize a cost functional, given initial and  final conditions where the controlled dynamics is given by nonholonomic mechanical system.  
In our paper, the controlled equations are derived using a basis
of vector fields adapted to  the nonholonomic distribution and the Riemannian
metric determined by the kinetic energy. 
Given a cost function, the
optimal control problem is understood as a 
\textit{constrained problem} or equivalently, under some mild regularity
conditions, as a Hamiltonian problem on the cotangent bundle of the nonholonomic distribution. A suitable
Lagrangian submanifold is also shown to lead to the correct 
dynamics. We demonstrate our techniques in several examples including 
a continuously variable transmission problem and motion
planning for obstacle avoidance problems.
\end{abstract}

\maketitle
\begin{center}
Dedicated to  H\'elene Frankowska and H\'ector J. Sussmann
\end{center}

\section{Introduction}

Although nonholonomic systems have been studied since the dawn of analytical mechanics, there has been some confusion over the correct formulation
of the equations of motion (see e.g. \cite{Bl}, \cite{BlMaZe}  and \cite{Manolo} for some 
of the history). Further it is only recently that their geometric
formulation has been understood. In addition, there has been recent interest
in the analysis of control problems for such systems.
 Nonholonomic control systems  exhibit distinctive features. In particular, 
many naturally underactuated systems are controllable, the controllability 
arising from the nonintegrability of the constraints. 

Nonholonomic optimal control problems arise in many  engineering applications, for instance systems with wheels, such as cars and  bicycles, and systems 
with blades or skates. There are thus multiple applications in the context of wheeled motion, space or mobile robotics and robotic manipulation.
In this paper, we will introduce  some new geometric techniques  in nonholonomic mechanics to study the case of force minimizing optimal control problems.


 The application of modern tools from differential
geometry in the fields of mechanics,  control theory, field theory and numerical integration  has led to 
significant progress in these research areas. For instance, the study
of the geometrical formulation of the nonholonomic equations of
motion has led to  better understanding of different engineering problems such  locomotion generation,
controllability, motion planning, and trajectory tracking (see e.g. \cite{Bl}, \cite{bloch2},
\cite{blochcrouch}, \cite{blochcrouch2}, \cite{blochcrouch3},
\cite{borisovmamaev}, \cite{bullolewis}, \cite{kellymurray},
\cite{jair},  \cite{LR1}, \cite{LMdD}, \cite{lewisandrew}, \cite{murraysastry},
\cite{ostrowoski} and references therein). Geometric techniques can also be used to study optimal control problems (see
\cite{blochcrouch3}, \cite{cortes noho}, \cite{CoMa}, \cite{Vu}, \cite{Vu2}, \cite{Su}, \cite{SuJu}).
Combining these ideas in this paper, we study  the underlying
geometry of optimal control problems for mechanical
systems subject to nonholonomic constraints and we apply it to several interesting examples.


 Classical nonholonomic constraints which are linear 
in the velocities can be geometrically encoded by a constant rank  distribution $\mathcal{D}$. As we will see, the  distribution $\mathcal{D}$ will play the role  of the velocity  phase space. Given a mechanical Lagrangian $L=K-V: TQ\rightarrow {\mathbb R}$ where $K$ and $V$ are the kinetic and potential energy, respectively.  and the distribution ${\mathcal D}$, the dynamics of the nonholonomic system is completely determined using the Lagrange-d'Alembert principle \cite{Bl}.
In this paper we will formulate a description in terms of a Levi Civita connection defined on the space of vector fields taking values on ${\mathcal D}$. This connection is obtained by projecting the standard Lie bracket using the Riemannian metric associated with the kinetic energy $K$ (see \cite{maria2}) and the typical characterization of the Levi-Civita connection (see also \cite{BlochHussein}). By adding controls in this setting we can study optimal control problems
 such as the force minimizing  problem.  Moreover, we can see that the dynamics of the optimal control problem is completely described by a Lagrangian submanifold of an appropriate cotangent bundle and, under some regularity conditions, the equations of motion are derived as classical  Hamilton's equations on the cotangent bundle of the distribution, $T^*{\mathcal D}$.
Although our approach is intrinsic, we also give a local description since it is important for working out examples. For this, it is necessary to choose
an adapted basis of vector fields for the distribution. From this point of view, we combine the techniques used previously by the authors of the paper (see  
\cite{maria2}, \cite{BloZen}, \cite{Blochmuru}). An additional advantage
 of our method is that symmetries may be naturally analyzed in this setting.

Concretely, the main results of our  paper can be summarized as follows:
\begin{itemize}
\item[$\bullet$] Geometric derivation of the equations of motion of nonholonomic  optimal control problems as a constrained problem on the tangent space to the constraint distribution ${\mathcal D}$. 

\item[$\bullet$] Construction of a Lagrangian submanifold representing the
dynamics of the optimal control problem and the corresponding Hamiltonian
representation when the system is regular.

\item[$\bullet$] Definition of a  Legendre transformation  establishing the
relationship  and correspondence between the Lagrangian and Hamiltonian
dynamics.

\item[$\bullet$] Application of our techniques to different examples
including optimal control of  the Chaplygin sleigh, 
a continuously variable transmission and motion planning for
obstacle avoidance problems.

\end{itemize}

\section{Nonholonomic mechanical systems}

Constraints on mechanical systems are  typically 
divided into two types: \textit{holonomic} and \textit{nonholonomic},
depending on whether the constraint can be derived from a constraint in
the configuration space or not. Therefore, the dimension of the
space of configurations is reduced by  holonomic constraints but not
by nonholonomic constraints. Thus, holonomic constraints allow a
reduction in the number of coordinates of the configuration space
needed to formulate a given problem (see~\cite{NF}).

We will restrict ourselves to the case  of nonholonomic constraints. Additionally, assume that  the constraints are given by a nonintegrable
distribution $\mathcal{D}$ on 
the configuration space $Q$. Locally, if we choose  local coordinates
$(q^i)$, $1\leq i\leq n=\dim Q$, the linear constraints on the
velocities are locally given by equations of the form
\[
\phi^{a}(q^i, \dot{q}^i)=\mu^a_i(q)\dot{q}^i=0, \quad 1\leq a\leq
m\leq n,
\]
depending, in general, on configuration coordinates and  their
velocities. {}From an intrinsic point of view, the linear
constraints are defined by a distribution ${\mathcal D}$ on $Q$ of
constant rank $n-m$ such that the annihilator of ${\mathcal D}$ is
locally given by
\[
{\mathcal D}^o = \operatorname{span}\left\{ \mu^{a}=\mu_i^{a}dq^i \; ; 1 \leq a
\leq m \right\}
\]
where the 1-forms $\mu^{a}$ are independent.

 In addition to these constraints, we
need to specify the dynamical evolution of the system, usually by
fixing a Lagrangian function $L\colon  TQ \to \R$. In mechanics, the
central concepts permitting the extension of mechanics from the
Newtonian point of view to the Lagrangian one are the notions of
virtual displacements and virtual work; these concepts were
originally formulated in the developments of mechanics in their application to
statics. In nonholonomic dynamics,  the procedure is given by the
\textit{Lagrange--d'Alembert principle}.
 This principle allows us to determine the set of possible values of the constraint forces from the set $\mathcal{D}$ of admissible kinematic states alone. The resulting equations of motion are
\begin{equation*}
\left[ \frac{d}{dt}\left( \frac{\partial L}{\partial \dot
q^i}\right) - \frac{\partial L}{\partial q^i} \right] \delta
q^i=0,
\end{equation*}
where $\delta q^i$ denotes the virtual displacements verifying
\begin{equation*}
\mu^a_i\delta q^i =0
\end{equation*}
(for the sake of simplicity, we will assume that the system is not
subject to non-conservative forces). This must be supplemented by
the constraint equations. By using the Lagrange multiplier rule, we
obtain
\begin{equation*}
\frac{d}{dt}\left( \frac{\partial L}{\partial \dot
q^i}\right)-\frac{\partial L}{\partial q^i}={\lambda}_a\mu^a_i  .
\end{equation*}
The term on the right hand side represents the  constraint force or
reaction force induced by the constraints. The functions $\lambda_a$
are Lagrange multipliers which, after being computed using the
constraint equations, allow us to obtain a set of second order
differential equations.

Now we restrict ourselves to the case of  nonholonomic mechanical
systems where the Lagrangian is of mechanical type
\[
L(v_q)=\frac{1}{2}\mathcal{G}(v_q, v_q) - V(q), \quad v_q\in T_qQ.
\]
Here $\mathcal{G}$ denotes a Riemannian metric on the configuration
space $Q$ representing the kinetic energy of the systems and
$V:Q\ra\R$ is a potential function. Locally, the metric is
determined by the matrix $M=(\mathcal{G}_{ij})_{1\leq i, j\leq n}$
where $\mathcal{G}_{ij}=\mathcal{G}(\partial/\partial q^i,
\partial/\partial q^j)$.

Denote by $\tau_{\mathcal{D}}:\mathcal{D}\ra Q$ the canonical
projection of $\mathcal{D}$ over $Q$ and
$\Gamma(\tau_{\mathcal{D}})$ the set of sections of $\tau_{D}$ which
 is just the set of vector fields $\mathfrak{X}(Q)$
taking values on $\mathcal{D}.$ If $X,Y\in\mathfrak{X}(Q),$ then
$[X,Y]$ denotes the standard Lie bracket of vector fields.

\begin{definition}\label{nonholonomicsystem}
A \textbf{nonholonomic mechanical system} on a manifold $Q$ is given
by the triple $(\mathcal{G}, V, \mathcal{D})$ where $\mathcal{G}$ is
a Riemannian metric on $Q,$ specifying the kinetic energy of the
system, $V:Q\ra\R$ is a smooth function representing the potential
energy and $\mathcal{D}$ a non-integrable distribution on $Q$
representing the nonholonomic constraints.
\end{definition}

\begin{remark}
Given $X,Y\in\Gamma(\tau_{\mathcal{D}})$ that is,
$X(x)\in\mathcal{D}_{x}$ and $Y(x)\in\mathcal{D}_{x}$ for all $x\in
Q,$ then it may happen that $[X,Y]\notin\Gamma(\tau_{\mathcal{D}})$
since $\mathcal{D}$ is nonintegrable. 
\end{remark}

We want to obtain a bracket defined for sections of $\mathcal{D}.$
Using the Riemannian metric $\mathcal{G}$ we can construct two
complementary orthogonal projectors
\begin{align*}
{\mathcal P}\colon &TQ\to {\mathcal D}\\
{\mathcal Q}\colon &TQ\to {\mathcal
D}^{\perp},
\end{align*} with respect to the tangent bundle orthogonal decomposition $\mathcal{D}\oplus\mathcal{D}^{\perp}=TQ.$

Therefore, given $X,Y\in\Gamma(\tau_{\mathcal{D}})$ we define the
\textbf{nonholonomic bracket}
$\lcf\cdot,\cdot\rcf:\Gamma(\tau_{\mathcal{D}})\times\Gamma(\tau_{\mathcal{D}})\rightarrow\Gamma(\tau_{\mathcal{D}})$
as
$$\lcf X, Y\rcf:=\mathcal{P}[X, Y]\; ,\qquad X, Y\in \Gamma(\tau_{\mathcal D})$$ (see
\cite{Balseiro},\cite{maria2},\cite{januzyleones}). It is clear that
this Lie bracket verifies the usual properties of a Lie bracket
except the Jacobi identity.
\begin{remark}
{\rm 
From a more differential geometric point of view, ${\mathcal D}$ with this modified bracket of sections inherits a skew-symmetric Lie algebroid structure  \cite{GrGrUr,maria2} where now the bracket of sections of the vector bundle does not satisfy in general the  Jacobi identity, as an expression of the nonintegrability of the distribution ${\mathcal D}$. 
}
\end{remark}

\begin{definition}
Consider the restriction of the Riemannian metric $\mathcal{G}$ to
the distribution $\mathcal{D}$
$$\mathcal{G}^{\mathcal{D}}:\mathcal{D}\times_{Q}\mathcal{D}\ra\R$$
and define the \textbf{Levi-Civita connection}
$$\nabla^{\mathcal{G}^{\mathcal{D}}}:\Gamma(\tau_{\mathcal{D}})\times\Gamma(\tau_{\mathcal{D}})\ra\Gamma(\tau_{\mathcal{D}})$$
determined by the following two properties:

\begin{enumerate}
\item $\lcf X,Y\rcf=\nabla_{X}^{\mathcal{G}^{\mathcal{D}}}Y-\nabla_{Y}^{\mathcal{G}^{\mathcal{D}}}X$\quad(Symmetry),\\
\item $X(\mathcal{G}^{\mathcal{D}}(Y,Z))=\mathcal{G}^{\mathcal{D}}(\nabla_{X}^{\mathcal{G}^{\mathcal{D}}}Y,Z)+\mathcal{G}^{\mathcal{D}}(Y,\nabla_{X}^{\mathcal{G}^{\mathcal{D}}}Z)$\quad
(Metricity).
\end{enumerate}

\end{definition}

Let $(q^{i})$ be coordinates on $Q$ and $\{e_{A}\}$ vector fields on
$\Gamma(\tau_{D})$ (that is, $e_A(x)\in {\mathcal D}_x$) such that $$\mathcal{D}_{x}=\hbox{span }\{e_{A}(x)\},
\quad x\in U\subset Q.$$ Then, we determine the \textbf{Christoffel
symbols} $\Gamma_{BC}^{A}$ of the connection
$\nabla^{\mathcal{G}^{\mathcal{D}}}$ by
$$\nabla_{e_{B}}^{\mathcal{G}^{\mathcal{D}}}e_{C}=\Gamma_{BC}^{A}(q)e_{A}.$$

\begin{definition}
A curve $\gamma:I\subset\R\ra\mathcal{D}$ is \textbf{admissible} if
there exists a curve $\sigma:I\subset\R\rightarrow Q$ projecting
$\gamma$ over $Q,$ that is, $\tau_{\mathcal{D}}\circ\gamma=\sigma;$
such that
$$\gamma(t)=\frac{d\sigma}{dt}(t).$$
\end{definition}

Given local coordinates on $Q,$ $(q^{i})$ $i=1,\ldots,n;$ and
$\{e_{A}\}$ a basis of sections on $\Gamma(\tau_{\mathcal{D}})$ such that
$\displaystyle{e_{A}=\rho_{A}^{i}(q)\frac{\partial}{\partial
q^{i}}}$ we introduce induced coordinates $(q^{i},y^{A})$ on
$\mathcal{D}$ where, if $e\in\mathcal{D}_{x}$ then
$e=y^{A}e_{A}(x).$ Therefore, $\gamma(t)=(q^{i}(t),y^{A}(t))$ is
admissible if
$$\dot{q}^{i}(t)=\rho_{A}^{i}(q(t))y^{A}(t).$$

Consider the restricted Lagrangian function
$\ell:\mathcal{D}\rightarrow\mathbb{R},$
$$\ell(v)=\frac{1}{2}\mathcal{G}^{\mathcal{D}}(v,v)-V(\tau_{D}(v)),\hbox{ with }  v\in\mathcal{D}.$$

\begin{definition}[\cite{maria2}]
A \textbf{solution of the nonholonomic problem} is an admissible
curve $\gamma:I\rightarrow\mathcal{D}$ such that
$$\nabla_{\gamma(t)}^{\mathcal{G}^{\mathcal{D}}}\gamma(t)+grad_{\mathcal{G}^{\mathcal{D}}}V(\tau_{\mathcal{D}}(\gamma(t)))=0.$$
\end{definition} Here the section $grad_{{\mathcal G}^{\mathcal{D}}}V\in\Gamma(\tau_{\mathcal{D}})$ is characterized by \[ {{\mathcal G}^{\mathcal{D}}}(grad_{{\mathcal
G}^{\mathcal{D}}}V, X) = X(V), \; \; \mbox{ for  every } X \in
\Gamma(\tau_{\mathcal{D}}).
\]

These equations are equivalent to the \textbf{nonholonomic
equations}. Locally, these are given by
\begin{eqnarray*}
\dot{q}^{i}&=&\rho_{A}^{i}(q)y^{A}\\
\dot{y}^{C}&=&-\Gamma_{AB}^{C}y^{A}y^{B}-(\mathcal{G}^{\mathcal{D}})^{CB}\rho_{B}^{i}\frac{\partial V}{\partial q^{i}}
\end{eqnarray*} where $(\mathcal{G}^{\mathcal{D}})^{AB}$ denotes the coefficients of the inverse matrix of $(\mathcal{G}^{\mathcal{D}})_{AB}$ where $\mathcal{G}^{\mathcal{D}}(e_{A},e_{B})=(\mathcal{G}^{\mathcal{D}})_{AB}.$

\begin{remark}
Observe that these equations only depend on the coordinates
$(q^{i},y^{A})$ on $\mathcal{D}$.  Therefore the nonholonomic
equations are free of Lagrange multipliers. These equations are
equivalent to the  \textit{nonholonomic Hamel equations} (see
\cite{BloZen}, \cite{Blochmuru} for example, and reference therein).

\end{remark}

\section{Optimal control of nonholonomic mechanical
systems}\label{section3}

The purpose of this section is to study optimal control problems for
a nonholonomic mechanical systems. We shall assume that all the
considered control systems are controllable, that is, for any two
points $q_0$ and $q_f$ in the configuration space $Q$, there exists
an admissible control $u(t)$ defined on the control manifold
$U\subseteq\R^{n}$ such that the system with initial condition $q_0$
reaches the point $q_f$ at time $T$ (see \cite{Bl,bullolewis} for
more details).

We will analyze the case when the dimension of the input or control distribution is equal to the rank of ${\mathcal D}$. If the rank of $\mathcal{D}$ is equal to the dimension of the control distribution, the system will be called a \textit{fully actuated nonholonomic system}.

\begin{definition}
A \textbf{solution of a fully actuated nonholonomic problem} is an
admissible curve $\gamma:I\rightarrow\mathcal{D}$ such that
$$\nabla_{\gamma(t)}^{\mathcal{G}^{\mathcal{D}}}\gamma(t)+grad_{\mathcal{G}^{\mathcal{D}}}V(\tau_{\mathcal{D}}(\gamma(t)))\in \Gamma (\tau_D),$$
or, equivalently, 
$$\nabla_{\gamma(t)}^{\mathcal{G}^{\mathcal{D}}}\gamma(t)+grad_{\mathcal{G}^{\mathcal{D}}}V(\tau_{\mathcal{D}}(\gamma(t)))=u^{A}(t)e_{A}(\tau_{\mathcal{D}}(\gamma(t)),$$
where $u^{A}$ are the control inputs.
\end{definition}

Locally, the equations may be  written as
\begin{eqnarray*}
\dot{q}^{i}&=&\rho_{A}^{i}y^{A}\\
\dot{y}^{C}&=&-\Gamma_{AB}^{C}y^{A}y^{B}-(\mathcal{G}^{\mathcal{D}})^{CB}\rho_{B}^{i}\frac{\partial V}{\partial q^{i}}+u^{C}.
\end{eqnarray*}

Given a cost function
\begin{eqnarray*}
C&:&\mathcal{D}\times U\rightarrow\mathbb{R}\\
&&(q^{i},y^{A}, u^A)\mapsto C(q^{i},y^{A},u^{A})
\end{eqnarray*} the \textit{optimal control problem} consists of finding an admissible curve
$\gamma:I\rightarrow\mathcal{D}$ solution of the fully actuated
nonholonomic problem given initial and final boundary conditions on
$\mathcal{D}$ and minimizing the functional
$$\mathcal{J}(\gamma(t),u(t)):=\int_{0}^{T}C(\gamma(t),u(t))dt,$$
where $\gamma$ is an admissible curve.

We define the submanifold $\mathcal{D}^{(2)}$ of $T\mathcal{D}$ by
\begin{equation}\label{D2}
\mathcal{D}^{(2)}:=\{v\in T\mathcal{D}\mid v=\dot{\gamma}(0)\hbox{
where } \gamma:I\rightarrow\mathcal{D} \hbox{ is admissible}\},
\end{equation} and we can choose coordinates $(x^{i},y^{A},\dot{y}^{A})$ on
$\mathcal{D}^{(2)}$ where the inclusion on $T\mathcal{D}$,
$i_{\mathcal{D}^{(2)}}:\mathcal{D}^{(2)}\hookrightarrow
T\mathcal{D}$ is given by
$$i_{\mathcal{D}^{(2)}}(q^{i},y^{A},\dot{y}^{A})=(q^{i},y^{A},\rho_{A}^{i}(q)y^{A},\dot{y}^{A}).$$ Therefore, $\mathcal{D}^{(2)}$ is locally described by the constraints on $T\mathcal{D}$ $$\dot{q}^{i}-\rho_{A}^{i}y^{A}=0.$$

Observe now that our optimal control problem is alternatively
determined by a smooth function
$\mathcal{L}:\mathcal{D}^{(2)}\rightarrow\mathbb{R}$ where

\begin{equation}\label{lagrangiancontrol}\mathcal{L}(q^{i},y^{A},\dot{y}^{C})=C\left(q^{i},
y^{A},\dot{y}^{C}+\Gamma_{AB}^{C}y^{A}y^{B}+(\mathcal{G}^{\mathcal{D}})^{CB}\rho_{B}^{i}\frac{\partial
V}{\partial q^{i}}\right).
\end{equation}

The following diagram summarizes the situation:

$$
\xymatrix{
\mathcal{D}^{(2)} \ar[dd]_{\mathcal{L}} \ar@{^{(}->}[rr]^{i_{{D}^{(2)}}} \ar[dr]^{\tau_{\mathcal{D}}^{(2,1)}} & \ & T\mathcal{D} \ar[dl]^{\tau_{T\mathcal{D}}} \ar[dr]^{T\tau_{\mathcal{D}}} \\
\ & \mathcal{D} \ar@{^{(}->}[rr]^{j} \ar[dr]_{\tau_{\mathcal{D}}} & \ & TQ \ar[dl]^{\tau_{TQ}} & \ \\
\mathbb{R} & \ & Q& \ &
}
$$ Here $j:\mathcal{D}\ra TQ$ is the canonical inclusion from $\mathcal{D}$ to $TQ,$ $\tau_{\mathcal{D}}^{(2,1)}:\mathcal{D}^{(2)}\ra\mathcal{D}$
and $\tau_{T\mathcal{D}}:T\mathcal{D}\ra\mathcal{D}$ are the
projections locally given by
$\tau_{\mathcal{D}}^{(2,1)}(q^{i},y^{A},\dot{y}^{A})=(q^{i},y^{A})$
and
$\tau_{T\mathcal{D}}(q^{i},y^{A},v^{i},\dot{y}^{A})=(q^{i},y^{A}),$
respectively. Finally, $T\tau_{\mathcal{D}}:T\mathcal{D}\ra TQ$ is
locally described as follows 
$(q^{i},y^{A},\dot{q}^{i},\dot{y}^{A})\mapsto (q^{i},\dot{q}^{i})$.

To derive the equations of motion for $\mathcal{L}$ we can use
standard variational calculus for systems with constraints defining
the extended Lagrangian $\widetilde{\mathcal{L}},$
$$\widetilde{\mathcal{L}}=\mathcal{L}+\lambda_{i}(\dot{q}^{i}-\rho_{A}^{i}y^{A}).$$
Therefore the equations of motion are
\begin{eqnarray}\label{eqELmulti}
\frac{d}{dt}\left(\frac{\partial\widetilde{\mathcal{L}}}{\partial\dot{q}^{i}}\right)-\frac{\partial\widetilde{\mathcal{L}}}{\partial q^{i}}&=&\dot{\lambda}_{i}-\frac{\partial\mathcal{L}}{\partial q^{i}}+\lambda_{j}\frac{\partial\rho_{A}^{j}}{\partial q^{i}}y^{A}=0,\nonumber\\
\frac{d}{dt}\left(\frac{\partial\widetilde{\mathcal{L}}}{\partial\dot{y}^{A}}\right)-\frac{\partial\widetilde{\mathcal{L}}}{\partial y^{A}}&=&\frac{d}{dt}\left(\frac{\partial \mathcal{L}}{\partial\dot{y}^{A}}\right)-\frac{\partial \mathcal{L}}{\partial y^{A}}+\rho_{A}^{i}\lambda_{i}=0,\\
\dot{q}^{i}&=&\rho_{A}^{i}y^{A}.\nonumber
\end{eqnarray}

\subsection{Example: continuously variable transmission
(CVT)}\label{cvt}

We want to study the optimal control of a simple model of a continuously
variable 
transmissions, where we assume that the belt cannot slip
(see \cite{Modin} for more details).   
\begin{center}
\begin{figure}[h]
  \includegraphics[width=8cm]{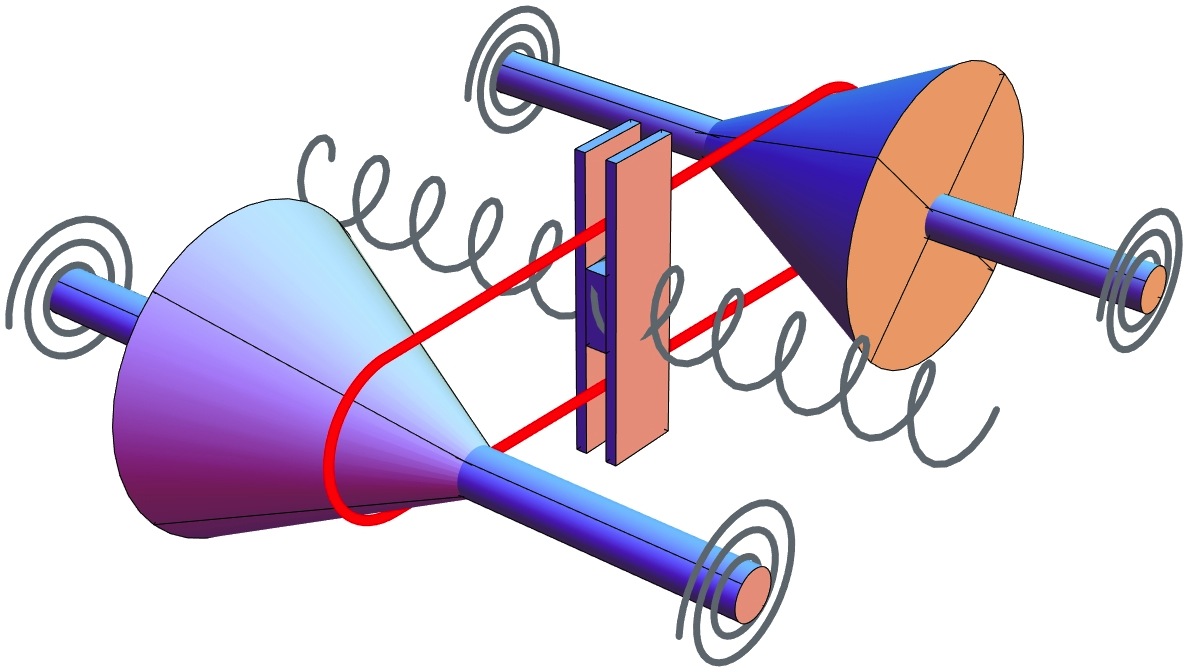}\
  \caption{Illustration of a continuously variable transmission  \cite{Modin}.}\end{figure}
\end{center}

The shafts are attached to spiral springs that are fixed to a chasis. The belt between the two cones is translated along the shafts in accordance with the coordinate $x$, thus providing a varying transmission ratio. The belt is kept in a plane perpendicular to the shafts, so that the belt keeps a constant length (see \cite{Modin} for a complete description and integrability of this system). The variables $\theta_1$ and $\theta_2$ denote the angular deflections of the shafts. $m$ denotes the mass of the belt slider, $J_1>0$ is the inertia
about the center of mass of the driving pulley and $J_2>0$ is the
inertia about the center of mass of the driven pulley. The
configuration space is $\mathbb{S}^1\times \mathbb{S}^1\times \R$
and the configuration is given by $q=(\theta_1,\theta_2,x)\in
\mathbb{S}^1\times \mathbb{S}^1\times \R$. 


The control inputs are denoted by $u_1$ and $u_2.$ The first one
corresponds to a force applied perpendicular to the center of mass
of the belt slider and the second one is the torque applied about
the the center of mass of the driving pulley.  Also, we assume that $x<1$ (which correspond to assuming that the gear ratio is finite). 

The belt imposes a constraint given by the no slip condition and is
expressed in differential form by
$$\omega=x\,d\theta_1-(1-x)\,d\theta_2.$$  Therefore the constraint
distribution $\mathcal{D}$ is given by
$$\mathcal{D}=\left\{\frac{1}{m}\frac{\partial}{\partial x},(1-x)\frac{\partial}{\partial \theta_1}+x\frac{\partial}{\partial \theta_2}\right\}.$$


The Lagrangian is metric on $\mathbb{S}^1\times \mathbb{S}^1\times
\R$ where the matrix associated with the metric $\mathcal{G}$ is
$$\mathcal{G}=\left(
                                      \begin{array}{ccc}
                                        J_1 & 0 & 0 \\
                                        0 & J_2 & 0 \\
                                        0 & 0 & m \\
                                      \end{array}
                                    \right).
$$ Then the Lagrangian $L:T(\mathbb{S}^1\times \mathbb{S}^1\times
\R)\ra\R$ is given by
$$L(q,\dot{q})=(\frac{J_1}{2}\dot{\theta}_1^{2}+\frac{J_2}{2}\dot{\theta}_2^{2})+\frac{m}{2}\dot{x}^{2}.$$

The projection map $\mathcal{P}:T(\mathbb{S}^1\times
\mathbb{S}^1\times \R)\ra\mathcal{D}$ is

\begin{eqnarray}\nonumber
&&\mathcal{P}(q,\dot{q})=\frac{J_1(1-x)^2}{J_1-2J_1x+J_1x^2+J_2x^2}d\theta_1\otimes\frac{\partial}{\partial \theta_1}+\frac{J_1x(1-x)}{J_1-2J_1x+J_1x^2+J_2x^2}d\theta_1\otimes\frac{\partial}{\partial \theta_2}\\
&&+\frac{J_2x(1-x)}{J_1-2J_1x+J_1x^2+J_2x^2}d\theta_2\otimes\frac{\partial}{\partial \theta_1}+\frac{J_2x^2}{J_1-2J_1x+J_1x^2+J_2x^2}d\theta_2\otimes\frac{\partial}{\partial \theta_2}+dx\otimes\frac{\partial}{\partial x}. \nonumber
\end{eqnarray}

Let $q=(\theta_1,\theta_2,x)$ be coordinates on the base manifold
$\mathbb{S}^1\times \mathbb{S}^1\times \R$ and take the basis
$\{X_{1},X_{2}\}$ of vector fields on $\mathbb{S}^1\times
\mathbb{S}^1\times \R.$ This basis induces adapted coordinates
$(\theta_1,\theta_2,x,y_{1},y_2)\in\mathcal{D}$  in the following way: Given the  vector fields $X_1$ and $X_2$ generating the distribution  $\mathcal{D}$ we obtain the relations for $q\in \mathbb{S}^1\times\mathbb{S}^1\times\mathbb{R}$
\begin{eqnarray*}
X_1(q)&=&\rho_1^1(q)\frac{\partial}{\partial \theta_1}+\rho_1^2(q)\frac{\partial}{\partial \theta_2}+\rho_1^3(q)\frac{\partial}{\partial x},\\
X_2(q)&=&\rho_2^1(q)\frac{\partial}{\partial \theta_1}+\rho_2^2(q)\frac{\partial}{\partial \theta_2}+\rho_2^3(q)\frac{\partial}{\partial x}.
\end{eqnarray*} Then, $$\rho_{1}^{1}= \rho_{1}^{2}=\rho_{2}^{3}=0,\quad\rho_{1}^{3}=\frac{1}{m},\quad
\rho_{2}^{1}=1-x,\quad\rho_{2}^{2}=x.
$$ Each element $e\in {\mathcal D}_q$ is expressed as a linear combination of these vector fields:
\[
e=y_1 X_1(q)+ y_2 X_2(q), \quad q\in \mathbb{S}^1\times\mathbb{S}^1\times\mathbb{R}.
\]

Therefore, the vector subbundle $\tau_{\mathcal{D}}:{\mathcal D}\rightarrow \mathbb{S}^1\times\mathbb{S}^1\times\mathbb{R}$ is locally described by the coordinates $(\theta_1,\theta_2,x; y_1, y_2)$; the first three for the base and the last two, for the fibers. 
Observe that 
\[
e=y_1\left(\frac{1}{m}\frac{\partial}{\partial  x}\right)+ y_2\left(
(1-x)\frac{\partial}{\partial \theta_1}+x\frac{\partial}{\partial \theta_2}\right)
\]
and, in consequence, ${\mathcal D}$ is described by the conditions (admissibility conditions): 
$$
\dot{\theta}_1=(1-x)y_{2},\quad
\dot{\theta}_2=xy_2,\quad
\dot{x}=\frac{1}{m}y_{1}
$$
as  a vector subbundle of $TQ$ where $y_1$ and $y_2$ are the adapted velocities relative to the basis of $\mathcal{D}$ defined before.

The nonholonomic
bracket is given by
$\ab=\mathcal{P}([\cdot,\cdot]).$ Observe now, \begin{eqnarray*}\lcf
X_1,X_2\rcf=\mathcal{P}[X_1,X_2]&=&\mathcal{P}\left(-\frac{1}{m}\frac{\partial}{\partial
\theta_1}+\frac{1}{m}\frac{\partial}{\partial
\theta_2}\right)\\
&=&-\frac{1}{m}\frac{J_{1}(1-x)-J_{2}x}{J_{2}x^{2}+J_{1}(1-x)^{2}}\left((1-x)\frac{\partial}{\partial\theta_{1}}+x\frac{\partial}{\partial\theta_{2}}\right).
\end{eqnarray*}






The restricted Lagrangian function in these new adapted coordinates
is rewritten as
$$\ell(\theta_1,\theta_2,x,y_1,y_2)=\frac{y_{2}^{2}}{2}((1-x)^{2}J_{1}+J_{2}x^{2})+\frac{1}{2m}y_1^2.$$

The Euler-Lagrange equations, together with  the admissibility conditions,
for this Lagrangian are

\begin{eqnarray*} \frac{\dot{y}_{1}}{m}&=&0,\quad
\dot{y}_{2}B(x)-\frac{y_{1}y_{2}A(x)}{m}=0\\
\dot{\theta}_{1}&=&(1-x)y_{2},\quad\dot{\theta}_{2}=xy_2,\quad\dot{x}=\frac{1}{m}y_1,
\end{eqnarray*}where $A(x)=J_1(1-x)-J_2x$ and $B(x)=(1-x)^{2}J_{1}+J_{2}x^{2}$.

Now, we add controls in our picture. Therefore the controlled
Euler-Lagrange equations are now
\begin{eqnarray*}
u_{1}&=&\dot{y}_{2}B(x)-\frac{y_{1}y_{2}A(x)}{m},\\
u_{2}&=&\frac{\dot{y}_{1}}{m},
\end{eqnarray*} together with
\begin{equation*}
\dot{\theta}_{1}=(1-x)y_{2},\quad\dot{\theta}_{2}=xy_2,\quad\dot{x}=\frac{1}{m}y_1.
\end{equation*}

The optimal control problem consists of finding an admissible curve
satisfying the previous equations given boundary conditions on
$\mathcal{D}$ and minimizing the functional

$$\mathcal{J}(\theta_1,\theta_2,x,y_1,y_2,u_1,u_2)=\frac{1}{2}\int_{0}^{T}\left(u_{1}^{2}+u_{2}^{2}
\right)dt,$$for the cost function
$C:\mathcal{D}\times U\ra\R$ given by
$$C(\theta_1,\theta_2,x,y_1,y_2,u_1,u_2)=\frac{1}{2}(u_1^2+u_2^2).$$

This optimal control problem is equivalent to the constrained
optimization problem determined by the lagrangian 
$\mathcal{L}:\mathcal{D}^{(2)}\ra\R$ given by
\begin{eqnarray*}
\mathcal{L}(\theta_1,\theta_2,x,y_1,y_2,\dot{y}_1,\dot{y}_2)&=&\frac{1}{2}\left(\dot{y}_{2}B(x)-\frac{y_{1}y_{2}A(x)}{m}\right)^{2}+\frac{\dot{y}_{1}^2}{2m^2}.
\end{eqnarray*}

Here, $\mathcal{D}^{(2)}$ is a submanifold of the vector bundle
$T\mathcal{D}$ over $\mathcal{D}$ defined by

\begin{small}$$\mathcal{D}^{(2)}:=\left\{(\theta_1,\theta_2,x,y_1,y_2,\dot{\theta}_1,\dot{\theta}_2,\dot{x},\dot{y}_1,\dot{y}_2)\in
T\mathcal{D}\Big{|}\dot{x}-\frac{1}{m}y_1=0,\dot{\theta}_1-(1-x)y_2=0,\dot{\theta}_2-xy_2=0\right\},$$\end{small}
where the inclusion
$i_{\mathcal{D}^{(2)}}:\mathcal{D}^{(2)}\hookrightarrow T\mathcal{D},$ is given by the map
$$i_{\mathcal{D}^{(2)}}(\theta_1,\theta_2,x,y_1,y_2,\dot{y}_1,\dot{y}_2)=
\left(\theta_1,\theta_2,x,y_1,y_2,(1-x)y_2,xy_2,\frac{y_1}{m},\dot{y}_1,\dot{y}_2\right).$$

The equations of motion for the extended Lagrangian
\begin{eqnarray*}
&&\widetilde{\mathcal{L}}(\theta_1,\theta_2,x,y_1,y_2,\dot{\theta}_1,\dot{\theta}_2,\dot{x},\dot{y}_1,\dot{y}_2,\lambda)=\\&&\mathcal{L}+\lambda_{1}\left(\dot{\theta}_1-(1-x)y_{2}\right)+\lambda_{2}\left(\dot{\theta}_2-xy_{2}\right)+\lambda_{3}\left(\dot{x}-\frac{1}{m}y_{1}\right)
\end{eqnarray*}
are

\begin{eqnarray*}
\dot{\lambda}_{1}&=&0,\quad \dot{\lambda}_{2}=0,\\
\dot{\lambda}_{3}&=&y_{2}(\lambda_1-\lambda_2)+\left(\dot{y}_2B(x)-A(x)\frac{y_1y_2}{m}\right)\left(\frac{y_1y_2(J_1+J_2)}{m}-2\dot{y}_2A(x)\right),\\
\lambda_{3}&=&-\frac{\ddot{y}_1}{m}-Ay_2\left(\dot{y}_2B(x)-A(x)\frac{y_1y_2}{m}\right),\\
0&=&\lambda_1(1-x)+\lambda_2x-\frac{3}{m}y_1A(x)\left(\dot{y}_2B(x)-A(x)\frac{y_1y_2}{m}\right)\\
&+&B(x)\left(\ddot{y}_2B(x)+\dot{y}_2\frac{2y_1A(x)}{m}-\frac{1}{m}\left(A(\dot{y}_1y_2+\dot{y}_2y_1)-\frac{y_1^{2}(J_1+J_2)}{m}\right)\right)
\end{eqnarray*} with

\begin{equation*}
\dot{\theta}_1=(1-x)y_{2},\quad\dot{\theta}_2=xy_2,\quad\dot{x}=\frac{1}{m}y_1.
\end{equation*}

%
%

The resulting system of equations for the optimal control problem of
the continuously variable transmission is difficult to solve
explicitly and from this observation it is clear that it is
necessary to develop numerical methods preserving the geometric
structure for these mechanical control systems. The construction of
geometric numerical methods for this kind of optimal control
problem is a  future research  topic as we
remark in \hyperref[conclusiones]{Section 6}.



\subsection{Example: the Chaplygin sleigh }\label{chaplygin}

We want to study the optimal control of the so-called
\textit{Chaplygin sleigh}  (see \cite{Bl}) introduced and studied in
1911 by Chaplygin \cite{Cha}, \cite{NF} and more recently by A. Ruina \cite{ruina} (see also \cite{Fedorov} and
\cite{FeZe05}). The sleigh is a rigid body moving on a horizontal
plane supported at three points, two of which slide freely without
friction while the third is a knife edge which allows no motion
orthogonal to its direction as show in Figure $3$.

We assume that the sleigh cannot move sideways. The configuration
space of this dynamical system is the group of Euclidean motions of
the two-dimensional plane $\R^2, SE(2),$ which we parameterize with
coordinates $(x, y, \theta)$ since an element $A\in SE(2)$ is
represented by the matrix$$\left(
                                             \begin{array}{ccc}
                                               \cos\theta & -\sin\theta & x \\
                                               \sin\theta & \cos\theta & y \\
                                               0 & 0 & 1 \\
                                             \end{array}
                                           \right) \hbox{ with
                                           }x,y\in\R\hbox{ and }\theta\in\mathbb{S}^{1}.
$$   $\theta$ and $(x, y)$ are the angular orientation of the
sleigh and position of the contact point of the sleigh on the plane,
respectively. Let $m$ be the mass of the
sleigh and $JI+ma^{2}$ is the inertia about the contact point,  where $I$ is the moment of inertia
about the center of mass $C$ and $a$ is the distance from the center of mass to the knife edge. The configuration space will be identified with  $\R^2\times\mathbb{S}^1$ with coordinates $q=(x,y,\theta)\in \R^2\times \mathbb{S}^1$.

\begin{center}
\begin{figure}[h]
  \includegraphics[width=8cm]{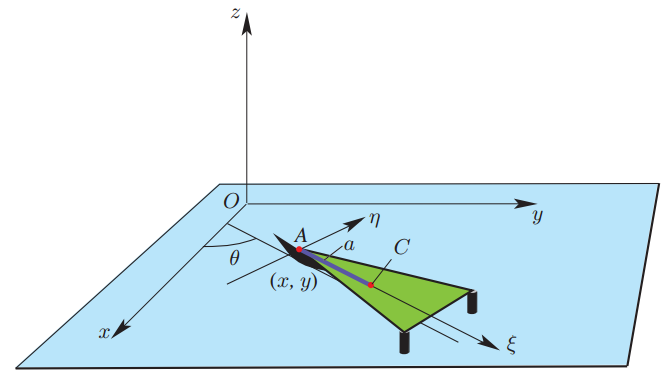}\label{sleigh}
  \caption{The Chaplygin sleigh}\end{figure}
\end{center}

The control inputs are denoted by $u_1$ and $u_2.$ The first one
corresponds to a force applied perpendicular to the center of mass
of the sleigh and the second one is the torque applied about the
vertical axis.

The constraint is given by the no slip condition and is
expressed in differential form by
$$\omega=\sin\theta\, dx-\cos\theta	, dy.$$  Therefore the constraint
distribution $\mathcal{D}$ is given by
$$\mathcal{D}=\left\{\frac{1}{J}\frac{\partial}{\partial\theta},\frac{\cos\theta}{m}\frac{\partial}{\partial x}+\frac{\sin\theta}{m}\frac{\partial}{\partial y}\right\}.$$

That is, the distribution is given by the span of the vector fields
\begin{eqnarray*}
X_{1}(q)&=&\frac{1}{J}\frac{\partial}{\partial\theta},\\
X_{2}(q)&=&\frac{\cos\theta}{m}\frac{\partial}{\partial x}+\frac{\sin\theta}{m}\frac{\partial}{\partial y}.
\end{eqnarray*}

The Lagrangian is metric on $Q$ where the matrix associated with the
metric $\mathcal{G}$ is  $$\mathcal{G}=\left(
                                      \begin{array}{ccc}
                                        m & 0 & 0 \\
                                        0 & m & 0 \\
                                        0 & 0 & J \\
                                      \end{array}
                                    \right).
$$ Then the Lagrangian $L:T(\R^2\times\mathbb{S}^{1})\ra\R$ is given by the kinetic energy of the body, which is a sum
of the kinetic energy of the center of mass and the kinetic energy
due to the rotation of the body $$L(q,\dot{q})=\frac{m}{2}(\dot{x}_C^{2}+\dot{y}_C^{2})+\frac{J}{2}\dot{\theta}^{2},$$ where $x_C = x+a\cos\theta$,$y_C = y+a\sin\theta$.

The projection map $\mathcal{P}:TQ\ra\mathcal{D}$ is
$$\mathcal{P}(q,\dot{q})=\cos^{2}\theta dx\otimes\frac{\partial}{\partial x}+\cos\theta\sin\theta dx\otimes\frac{\partial}{\partial y}+\cos\theta\sin\theta dy\otimes\frac{\partial}{\partial x}+\sin^{2}\theta dy\otimes\frac{\partial}{\partial y}+d\theta\otimes\frac{\partial}{\partial\theta}.$$


Let $q=(x,y,\theta)$ be coordinates on the base manifold $\R^2\times\mathbb{S}^1$ and
take the basis $\{X_{1},X_{2}\}$ of vector fields of $\mathcal{D}.$ This basis
induces adapted coordinates $(x,y,\theta,y_1,y_2)\in\mathcal{D}$ in the following way: Given the  vector fields $X_1$ and $X_2$ generating the distribution  we obtain the relations for $q\in \R^2\times\mathbb{S}^1$
\begin{eqnarray*}
X_1(q)&=&\rho_1^1(q)\frac{\partial}{\partial x}+\rho_1^2(q)\frac{\partial}{\partial y}+\rho_1^3(q)\frac{\partial}{\partial\theta},\\
X_2(q)&=&\rho_2^1(q)\frac{\partial}{\partial x}+\rho_2^2(q)\frac{\partial}{\partial y}+\rho_2^3(q)\frac{\partial}{\partial\theta}.
\end{eqnarray*} Then, $$\rho_{1}^{1}= \rho_{1}^{2}=\rho_{2}^{3}=0,\quad\rho_{1}^{3}=\frac{1}{J},\quad
\rho_{2}^{1}=\frac{\cos\theta}{m},\quad\rho_{2}^{2}=\frac{\sin\theta}{m}.
$$ Each element $e\in {\mathcal D}_q$ is expressed as a linear combination of these vector fields:
\[
e=y_1 X_1(q)+ y_2 X_2(q), \quad q\in \R^2\times\mathbb{S}^1.
\]

Therefore, the vector subbundle $\tau_{\mathcal{D}}:{\mathcal D}\rightarrow \R^2\times\mathbb{S}^1$ is locally described by the coordinates $(x,y,\theta; y_1, y_2)$; the first three for the base and the last two, for the fibers. 
Observe that 
\[
e=y_1\left(\frac{1}{J}\frac{\partial}{\partial \theta}\right)+ y_2\left(
\frac{\cos\theta}{m}\frac{\partial}{\partial x}+\frac{\sin\theta}{m}\frac{\partial}{\partial y}\right)
\]
and, in consequence, ${\mathcal D}$ is described by the conditions (admissibility conditions): 
$$
\dot{x}=\frac{\cos\theta}{m}y_2,\quad
\dot{y}=\frac{\sin\theta}{m}y_2,\quad
\dot{\theta}=\frac{1}{J}y_1
$$
as  a vector subbundle of $TQ$ where $y_1$ and $y_2$ are the adapted velocities relative to the basis of $\mathcal{D}$ defined before.

The nonholonomic bracket given
by $\lcf\cdot,\cdot\rcf=\mathcal{P}([\cdot,\cdot])$ satisfies $$\lcf
X_1,X_2\rcf=\mathcal{P}[X_1,X_2]=\mathcal{P}\left(-\frac{1}{Jm}\sin\theta\frac{\partial}{\partial
x}+\frac{\cos\theta}{Jm}\frac{\partial}{\partial y}\right)=0.$$

The restricted Lagrangian function in the new adapted coordinates is
given by
$$\ell(x,y,\theta,y_1,y_2)=\frac{1}{2m}(y_2)^2+\frac{b}{2J}(y_1)^2 \hbox{ where } b=\frac{a^2m}{J}.$$

Therefore, the equations of motion are
$$
\frac{b\dot{y}_1}{J}=0,\quad \frac{\dot{y}_2}{m}=0,\quad
\dot{x}=\frac{\cos\theta}{m}y_2,\quad\dot{y}=\frac{\sin\theta}{m}y_2,\quad\dot{\theta}=\frac{1}{J}y_1.
$$

Now, by adding controls in our picture, the controlled Euler-Lagrange
equations are written as
$$\frac{b\dot{y}_1}{J}=u_2,\quad
\frac{\dot{y}_2}{m}=u_1,\quad
\dot{x}=\frac{\cos\theta}{m}y_2,\quad\dot{y}=\frac{\sin\theta}{m}y_2,\quad\dot{\theta}=\frac{1}{J}y_1.
$$

The optimal control problem consists on finding an admissible curve
satisfying the previous equations given boundary conditions on
$\mathcal{D}$ and minimizing the functional
$\mathcal{J}(x,y,\theta,y_1,y_2,u_1,u_2)=\frac{1}{2}\int_{0}^{T}\left(u_{1}^{2}+u_{2}^{2}
\right)dt$, for the cost function $C:\mathcal{D}\times U\ra\R$
given by
\begin{equation}\label{costchaplygin}C(x,y,\theta,y_1,y_2,u_1,u_2)=\frac{1}{2}(u_1^2+u_2^2).\end{equation}

As before, the optimal control problem is equivalent to solving the
constrained  optimization problem determined by
$\mathcal{L}:\mathcal{D}^{(2)}\ra\R,$ where

$$
\mathcal{L}(x,y,\theta,y_1,y_2,\dot{y}_1,\dot{y}_2)=\frac{1}{2}\left(\frac{b^2\dot{y}_1^2}{J^2}+\frac{\dot{y}_2^2}{m^2}\right).$$ Here, $\mathcal{D}^{(2)}$ is a submanifold of the vector bundle
$T\mathcal{D}$ over $\mathcal{D}$ defined by

$$\mathcal{D}^{(2)}:=\left\{(x,y,\theta,y_1,y_2,\dot{x},\dot{y},\dot{\theta},\dot{y}_1,\dot{y}_2)\in T\mathcal{D}\Big{|}\dot{x}-\frac{\cos\theta}{m}y_2=0,\dot{y}-\frac{\sin\theta}{m}y_2=0,\dot{\theta}-\frac{1}{J}y_1=0\right\},$$
where the inclusion
$i_{\mathcal{D}^{(2)}}:\mathcal{D}^{(2)}\hookrightarrow
T\mathcal{D},$ is given by the map
$$i_{\mathcal{D}^{(2)}}(x,y,\theta,y_1,y_2,\dot{y}_1,\dot{y}_2)=
\left(x,y,\theta,y_1,y_2,\frac{\cos\theta}{m}y_2,\frac{\sin\theta}{m}y_2,\frac{1}{J}y_1,\dot{y}_1,\dot{y}_2\right).$$

The equations of motion for the  extended Lagrangian
\begin{eqnarray*}\widetilde{\mathcal{L}}(x,y,\theta,y_1,y_2,\dot{x},\dot{y},\dot{\theta},\dot{y}_1,\dot{y}_2,\lambda)&=&\mathcal{L}+\lambda_{1}\left(\dot{x}-\frac{\cos\theta}{m}y_{2}\right)+\lambda_{2}\left(\dot{y}-\frac{\sin\theta}{m}y_{2}\right)\\
&&+\lambda_{3}\left(\dot{\theta}-\frac{1}{J}y_{1}\right)
\end{eqnarray*}
are
\begin{eqnarray*}
\dot{\lambda}_{1}&=&0,\quad \dot{\lambda}_{2}=0,\quad
\dot{\lambda}_{3}=\frac{y_{2}}{m}\left(\lambda_{1}\sin\theta-\lambda_{2}\cos\theta\right),\\
\lambda_{3}&=&-\frac{b^2\ddot{y}_{1}}{J},\quad
\ddot{y}_{2}=-m(\lambda_{1}\cos\theta+\lambda_2\sin\theta)
\end{eqnarray*} with

\begin{equation*}
\dot{x}=\frac{\cos\theta}{m}y_2,\quad\dot{y}=\frac{\sin\theta}{m}y_2,\quad\dot{\theta}=\frac{1}{J}y_1.
\end{equation*}

The first two equations can be integrated as $\lambda_{1}=c_1$ and
$\lambda_{2}=c_2$ where $c_1$ and ·$c_2$ are constants and
differentiating the equation for $\lambda_3$ with respect to the time and
substituting into the third equation, the problem is reduced to solve the
system
$$
\frac{\dddot{y}_{1}}{J}=\frac{y_{2}}{mb^{2}}\left(c_{2}\cos\theta-c_{1}\sin\theta\right),\quad
\ddot{y}_{2}=-m(c_{1}\cos\theta+c_2\sin\theta),
$$ with
\begin{equation*}
\dot{x}=\frac{\cos\theta}{m}y_2,\quad\dot{y}=\frac{\sin\theta}{m}y_2,\quad\dot{\theta}=\frac{1}{J}y_1.
\end{equation*}

If we suppose, $\lambda_1=0,\lambda_2=0$ (that is, $c_1=c_2=0$) then
the system can be reduced to solve
$$
\dddot{y}_{1}=0\ \hbox{    and   }\
\ddot{y}_{2}=0.
$$
Integrating these equations and using the admissibility conditions
we obtain constants of integration  $c_i,$ $i=3,\ldots,8$ and the 
equations
\begin{eqnarray*}
\theta(t)&=&\frac{c_3t^3}{6J}+\frac{c_4t^2}{2J}+\frac{c_5t+c_6}{J},\\
x(t)&=&\frac{1}{m}\int_{0}^{t}\cos\left(\frac{c_3s^3+3c_4s^2+6c_5s+6c_6}{6J}\right)\left(c_7s+c_8\right)ds,\\
y(t)&=&\frac{1}{m}\int_{0}^{t}\sin\left(\frac{c_3s^3+3c_4s^2+6c_5s+6c_6}{6J}\right)\left(c_7s+c_8\right)ds.
\end{eqnarray*}

Therefore the controls $u_1$ and $u_2$ are
$$u_1(t)=\frac{c_7}{m},\quad
u_2(t)=\frac{c_3t+c_4}{J}.$$

\begin{remark}\label{comparacion}
A similar optimal control problem was studied also \cite{BlochHussein}.
The authors have also used the theory of affine connections to analyze
the optimal control problem of underactuated nonholonomic mechanical
systems. The main difference with our approach  is that in our paper
 we are working on the distribution $\mathcal{D}$ itself. We impose the extra condition
$\lambda_1=\lambda_2=0$ to obtain  explicitlly the controls minimizing
the cost function.
Usually, there is prescribed an initial boundary condition
on $\mathcal{D}$ and a final boundary condition on $\mathcal{D}$. For
the Chaplygin sleigh we impose conditions
$(x(0),y(0),\theta(0),y_1(0),y_2(0))$ and
$(x(T),y(T),\theta(T),y_1(T),y_2(T))$. Heuristically, observe that
if we transform these conditions into initial conditions we will
need to take the initial condition\\
$(x(0),y(0),\theta(0),y_1(0),y_2(0),\dot{y}_1(0),\dot{y}_2(0),\lambda_1(0),\lambda_2(0),\lambda_{3}(0))$
and it is not necessary  that some of the multipliers are zero from the
very beginning. 
\end{remark}

\subsection{Application to motion planing for obstacle avoidance: The Chaplygin sleigh with obstacles}

In this section, we use the same model of the Chaplygin sleigh from
the previous section to show how obstacle avoidance can be achieved
with our approach using navigation functions. A
navigation function is a potential field-based function used to
model an obstacle as a repulsive area or surface
\cite{Khatib},\cite{Koditschek}.

For the Chaplygin sleigh, consider the following boundary
conditions on the distribution $\mathcal{D}$: $x(0)=0,\quad y(1)=0,\quad \theta(0)=0,\quad y_1(0)=0,\quad y_2(0)=0$ and $x(T)=1,\quad y(T)=1,\quad \theta(T)=0,\quad y_1(T)=0,\quad y_2(T)=0.$

Let the obstacle be circular in the $xy$-plane, located at the
point $(x_C,y_C)=(0.5,0.5)$. For llustrative purposes, we use a simple
inverse square law for the navigation function. Let $V(x,y)$ given
by $$V(x,y)=\frac{\kappa}{(x-x_C)^{2}+(y-y_C)^{2}}$$ where the
parameter $\kappa$ is introduced to control the strength of the
potential function. 

Appending the potential into the cost functional \eqref{costchaplygin}
the optimal control problem is equivalent to solve the constrained
optimization problem determined by
$\mathcal{L}:\mathcal{D}^{(2)}\ra\R,$ where

$$
\mathcal{L}(x,y,\theta,y_1,y_2,\dot{y}_1,\dot{y}_2)=\frac{b^2\dot{y}_1^2}{2J^2}+\frac{\dot{y}_2^2}{2m^2}+\frac{\kappa}{2((x-x_C)^{2}+(y-y_C)^2)}.$$

The equations of motion for the extended Lagrangian
\begin{eqnarray*}\widetilde{\mathcal{L}}(x,y,\theta,y_1,y_2,\dot{x},\dot{y},\dot{\theta},\dot{y}_1,\dot{y}_2,\lambda)&=&\mathcal{L}+\lambda_{1}\left(\dot{x}-\frac{\cos\theta}{m}y_{2}\right)+\lambda_{2}\left(\dot{y}-\frac{\sin\theta}{m}y_{2}\right)\\
&&+\lambda_{3}\left(\dot{\theta}-\frac{1}{J}y_{1}\right)
\end{eqnarray*}
are
\begin{eqnarray*}
\dot{\lambda}_{1}&=&-\frac{\kappa(x-x_C)}{((x-x_C)^{2}+(y-y_C)^{2})^{2}}\quad \dot{\lambda}_{2}=-\frac{\kappa(y-y_C)}{((x-x_C)^{2}+(y-y_C)^{2})^{2}},\\
\dot{\lambda}_{3}&=&\frac{y_{2}}{m}\left(\lambda_{1}\sin\theta-\lambda_{2}\cos\theta\right),\quad\lambda_{3}=-\frac{b^2\ddot{y}_{1}}{J},\quad
\ddot{y}_{2}=-m(\lambda_{1}\cos\theta+\lambda_2\sin\theta)
\end{eqnarray*} with

\begin{equation*}
\dot{x}=\frac{\cos\theta}{m}y_2,\quad\dot{y}=\frac{\sin\theta}{m}y_2,\quad\dot{\theta}=\frac{1}{J}y_1.
\end{equation*}

We solve the earlier boundary value problem for several values of $\kappa$.  Starting with $\kappa=0$, which corresponds to a zero potential function, we incremente $\kappa$ until the potential field was strong enough to prevent the sleigh from interfering with the obstacle. We try with $\kappa=0, 0.01, 0.1, 0.25,$ and $0.5$ for $T$=1. The result is shown in Fig. $4$. Note that for $\kappa=0.25$ and $0.5$ the sleigh avoids the obstacle, and as one may anticipate, as $\kappa$ increases, the total control effort and therefore, the total cost  $\displaystyle{\mathcal{J}=\frac{1}{2}\int_0^1(u_1^2+u_2^2+V(x,y))dt}$ increases.  For example, $\mathcal{J}=17.0242$ when $\kappa=0.25$ and $\mathcal{J}=18.4634$ when $\kappa=0.5$. Hence, we select $\kappa=0.25$ since it corresponds to a trajectory that avoids the obstacle with the least possible cost (of all five tried in this simulation). The trajectories profile is shown in Figures $5$, $6$ and $7$. This example illustrate how our approach can be used with the method of navigation functions of optimal motion generation for obstacle avoidance.

\begin{center}
\begin{figure}[h]
  \centering\includegraphics[width=9.5cm]{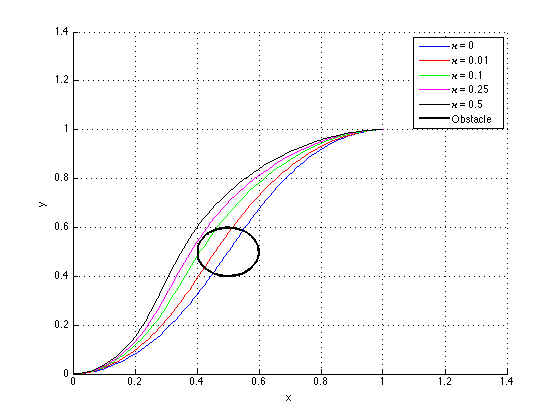}\label{fig4}
\caption{The extremals solving the boundary value problem with $\kappa=0,0.01,0.1,0.25$ and $0.5$.}\end{figure}
\end{center}

\begin{center}
\begin{figure}[h]
  \centering\includegraphics[width=7.5cm]{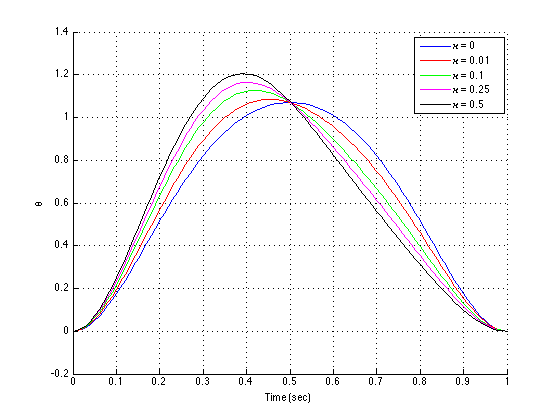}
\caption{Behavior of $\theta$ for $\kappa=0,0.01,0.1,0.25$ and $0.5$.}\end{figure}
\end{center}

\begin{figure}[h]
\includegraphics[width=7.3cm]{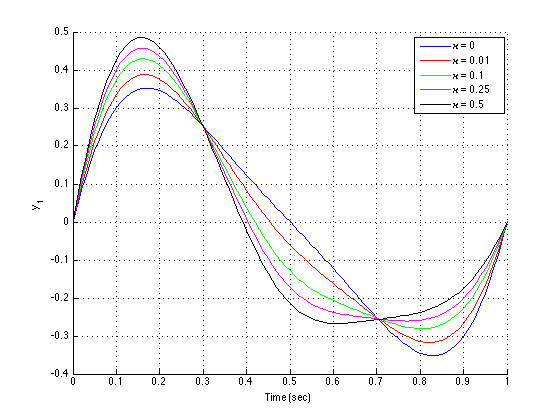}
\includegraphics[width=7.3cm]{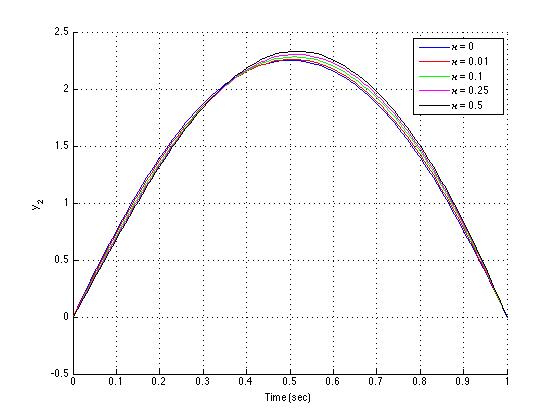}
\caption{Behavior of the velocites $y_1$ (left) and $y_2$ (right) for $\kappa=0,0.01,0.1,0.25$ and $0.5$}
\end{figure}

\begin{figure}[h]
\includegraphics[width=7.3cm]{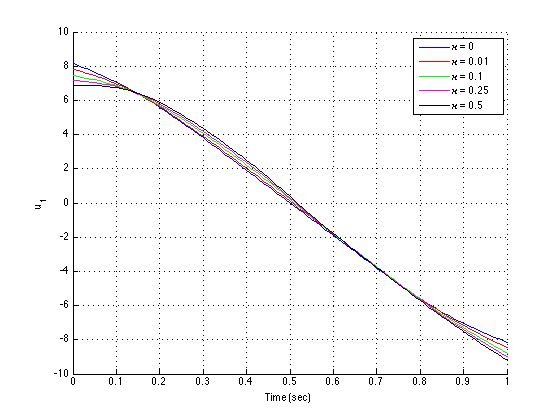}
\includegraphics[width=7.3cm]{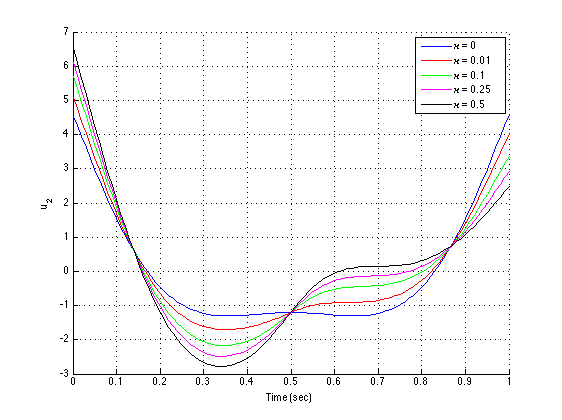}
\caption{Behavior of the controls $u_1$ (left) and $u_2$ (right) for $\kappa=0,0.01,0.1,0.25$ and $0.5$}
\end{figure}

\section{Lagrangian submanifolds and nonholonomic optimal control problems}
In this section we study the construction of Lagrangian submanifold representing intrinsically the dynamics of the optimal control problem and the corresponding Hamiltonian representation when the system is regular. In the regular case, the definition of a particular Legendre transformation give rise the relationship and correspondence between the Lagrangian and Hamiltonian dynamics.
\subsection{Lagrangian submanifolds}
In this subsection we will construct Lagrangian submanifolds that are interesting for our purposes in the study of the geometry of optimal control problems of controlled mechanical systems 
(see \cite{LiMa,Wein}).

\begin{definition}
Given a finite-dimensional symplectic manifold $(P, \omega)$ and a
submanifold $N$, with canonical inclusion $i_N: N\hookrightarrow P$,
 $N$ is said to be a \textbf{Lagrangian submanifold}  if
 $i^{*}_{N}\,\omega=0$ and $$\mbox{dim}\hspace{1mm}N=\frac{1}{2}\mbox{dim}\hspace{1mm}P.$$
\end{definition}

A distinguished symplectic manifold is the cotangent bundle  $T^*Q$
of any manifold $Q$. If we choose local  coordinates  $(q^i)$,
$1\leq i\leq n$, then $T^*Q$ has induced  coordinates $(q^i, p_i)$.
Denote by $\pi_Q: T^*Q\to Q$ the canonical projection of the
cotangent bundle defined by $\pi_Q(\epsilon_q)=q$, where
$\epsilon_q\in T^*_{q}Q$. Define the Liouville 1-form or canonical
1-form $\theta_Q\in \Lambda^1(T^*Q)$ by
\[
\bra(\theta_Q)_{\epsilon}\,,\,X\ket=\bra\epsilon\,,\,T\pi_Q(X)\ket, \hbox{  where  } X\in T_{\epsilon}T^*Q\; ,\ \epsilon\in T^*Q.
\]
In local coordinates we have that $\theta_Q=p_i\, \mathrm{d}q^i$.
The canonical two-form $\omega_Q$ on $T^*Q$ is the symplectic form
$\omega_Q=-\mathrm{d}\theta_Q$ (that is
$\omega_{Q}=\mathrm{d}q^{i}\wedge\mathrm{d}p_{i}$).

Now, we will introduce some special Lagrangian submanifolds of the
symplectic manifold $(T^*Q, \omega_Q)$. For instance, the image
${\Sigma}_{\lambda}=\lambda(Q)\subset T^*Q$ of a closed 1-form
$\lambda\in \Lambda^1Q$ is a Lagrangian submanifold of $(T^*Q,
\omega_Q)$, since $\lambda^* \omega_Q = - \mathrm{d}\lambda=0$. We
then obtain a submanifold diffeomorphic to $Q$ and transverse to the
fibers of $T^*Q$. When $\lambda$ is exact, that is,
$\lambda=\mathrm{d}f$, where $f: Q\to \R$, we say that $f$ is a
\emph{generating} function of the Lagrangian submanifold
${\Sigma}_{\lambda}=\Sigma_{f}$ (see \cite{Wein}).

A useful extension of the previous construction is the following
result due to W.Tulczyjew:

\begin{theorem}[\cite{Tulczy1},\cite{Tulczy2}]\label{tulchi}
  \label{thm:tulczyjew}
  Let $Q$ be a smooth manifold, $\tau_Q:TQ\Flder Q$ its tangent bundle projection, $ N \subset Q $ a submanifold, and $ f
  \colon N \rightarrow \mathbb{R} $.  Then
  \begin{multline*}
    \Sigma _f = \bigl\{ p \in T ^\ast Q \mid \pi _Q (p) \in N \text{
        and } \left\langle p, v \right\rangle = \left\langle
        \mathrm{d} f , v \right\rangle \\
      \text{ for all } v \in T N \subset T Q \text{ such that } \tau
      _Q (v) = \pi _Q (p) \bigr\}
  \end{multline*}
  is a Lagrangian submanifold of $ T ^\ast Q $.
\end{theorem}
Taking $f$ as the zero function, for example, we obtain the
following Lagrangian submanifold
\[
\Sigma_{0}=\lc p\in T^{*}Q\big|_{N}\,|\,\bra p\,,\,v\ket=0\,,\,\forall\,v\in TN\,\mbox{ with }\,\tau_{Q}(v)=\pi_{Q}(p)\rc,
\]
which is just the {\bf conormal bundle} of $N$:
\[
\nu^{*}(N)=\lc p\in T^{*}Q\big|_{N}\,;\,p\big|_{T_{\pi(p)}N}=0\rc.
\]

\subsection{Lagrangian submanifold description of nonholonomic
mechanical control problems} Next, we derive  the equations of motion representing the dynamics of the optimal control problem .

Given the function $\mathcal{L}:\mathcal{D}^{(2)}\ra\R, $ following
Theorem \eqref{tulchi}, when $N=\mathcal{D}^2\subset T{\mathcal D}$ we have the Lagrangian
submanifold
$\Sigma_{\mathcal{L}}\subset
T^{*}T\mathcal{D}.$ Therefore,
$\mathcal{L}:\mathcal{D}^{(2)}\ra\mathbb{R}$ generates a Lagrangian
submanifold $\Sigma_{\mathcal{L}}\subset T^{*}T\mathcal{D}$ of the
symplectic manifold $(T^{*}T\mathcal{D},\omega_{T\mathcal{D}})$
where $\omega_{T\mathcal{D}}$ is the canonical symplectic 2-form on
$T^{*}T\mathcal{D}.$

The relationship between these spaces is summarized in the following
diagram:
$$
\xymatrix{
\Sigma_{\mathcal{L}} \ar[d]_{(\pi_{T^{*}T\mathcal{D}})\mid_{\Sigma_{\mathcal{L}}}} \ar@{^{(}->}[rr]^{i_{\Sigma_L}} & \ & T^{*}T\mathcal{D} \ar[d]^{\pi_{T^{*}T\mathcal{D}}} \\
\mathcal{D}^{(2)}\ar@{^{(}->}[rr]^{i_{{D}^{(2)}}} & \ & T\mathcal{D}
}
$$

\begin{proposition}
Let $\mathcal{L}:\mathcal{D}^{(2)}\rightarrow\mathbb{R}$ be a 
 $C^{\infty}$-function. Consider the inclusion
$i_{\mathcal{D}^{(2)}}:\mathcal{D}^{(2)}\ra T\mathcal{D}$ where
$\omega_{T\mathcal{D}}$ is the canonical symplectic 2-form in
$T^{*}T\mathcal{D}.$ Then
$$\Sigma_\mathcal{L}=\{\mu\in T^{*}T\mathcal{D}|i_{\Sigma_{\mathcal{L}}}^{*}\mu=d\mathcal{L}\}\subset
T^{*}T\mathcal{D}$$ is a Lagrangian submanifold of
$(T^{*}T\mathcal{D},\omega_{T\mathcal{D}}).$
\end{proposition}

\begin{definition} Let $\mathcal{D}$ be a non-integrable distribution, $T\mathcal{D}$ its tangent bundle and $\mathcal{D}^{(2)}$ the subbundle of $T\mathcal{D}$ defined on \eqref{D2}.
A \textbf{second-order nonholonomic system} is a pair
$(\mathcal{D}^{(2)},\Sigma_{\mathcal{L}})$ where
$\Sigma_\mathcal{L}\subset T^{*}T\mathcal{D}$ is the Lagrangian
submanifold generated by $\mathcal{L}:\mathcal{D}^{(2)}\ra\R.$
\end{definition}

Consider local coordinates $(q^{i},y^{A},\dot{q}^{i},\dot{y}^{A})$
on $T\mathcal{D}.$ These coordinates induce local coordinates
$(q^{i},y^{A},\dot{q}^{i},\dot{y}^{A},\mu_{i},\mu_{A},\gamma_{i},\gamma_{A})$
on $T^{*}T\mathcal{D}.$ Therefore, locally, the system is
characterized by the following set of equations on
$T^{*}T\mathcal{D}$

\begin{eqnarray}\label{eqlagsub}
\mu_{i}+\gamma_{j}\frac{\partial\rho_{A}^{j}}{\partial q^{i}}y^{A}&=&\frac{\partial\mathcal{L}}{\partial q^{i}},\nonumber\\
\mu_{A}+\gamma_{j}\rho_{A}^{j}&=&\frac{\partial\mathcal{L}}{\partial y^{A}},\\
\gamma_{A}&=&\frac{\partial\mathcal{L}}{\partial\dot{y}^{A}},\nonumber\\
\dot{q}^{i}&=&\rho_{A}^{i}y^{A}.\nonumber
\end{eqnarray}

\begin{remark}
Typically local coordinates on $\Sigma_{\mathcal{L}}\subset
T^{*}T\mathcal{D}$ are $(q^{i},y^{A},\dot{y}^{A},\gamma_{i})$ where
$\gamma_{i}$ plays the role of Lagrange multipliers.
\end{remark}

\begin{remark}In the case of the Chaplygin sleigh local coordinates on $T^{*}T\mathcal{D}$ will
be given by
$(x,y,\theta,y_1,y_2,\dot{x},\dot{y},\dot{\theta},\dot{y}_1,\dot{y}_2,\mu_x,\mu_y,\mu_{\theta},\mu_1,\mu_2,\gamma_x,\gamma_y,\gamma_{\theta},\gamma_1,\gamma_2)$,
where
$(x,y,\theta,y_1,y_2,\dot{x},\dot{y},\dot{\theta},\dot{y}_1,\dot{y}_2)$
are local coordinates on $T\mathcal{D}$. The Lagrangian submanifold
of  $T^{*}T\mathcal{D}$ is described by the equations 
\begin{eqnarray*}
\mu_{x}&=&0,\quad\mu_y=0,\\
\mu_{\theta}&=&\frac{y_2}{m}\left(\gamma_x\sin\theta-\gamma_y\cos\theta\right),\\
\mu_1&=&-\frac{b^2\gamma_{\theta}}{J},\quad\mu_2=-m(\gamma_x\cos\theta+\gamma_y\cos\theta),\\
\gamma_1&=&\frac{b^2\dot{y}_{1}}{J^2},\quad\gamma_2=\frac{\dot{y}_2}{m^{2}},\\
\dot{x}&=&\frac{\cos\theta}{m}y_2,\quad\dot{y}=\frac{\sin\theta}{m}y_2,\quad\dot{\theta}=\frac{y_1}{J}.
\end{eqnarray*} After a straightforward computation one can check easily that these equations are equivalent with those obtained in the Lagrangian formalism.
\end{remark}

\subsection{Legendre transformation and regularity condition}

We define the map $\Psi:T^{*}T\mathcal{D}\ra T^{*}\mathcal{D}$ as
$$\langle\Psi(\mu_{v_{x}}),X(x)\rangle=\langle\mu_{v_{x}},X^{V}(v_{x})\rangle,$$ where
$\mu\in T^{*}T\mathcal{D},v_{x}\in T_{x}\mathcal{D},$ $X(x)\in
T_{x}\mathcal{D}$ and $X^{V}(v_x)\in T_{v_x}T\mathcal{D}$ is its
vertical lift to $v_x.$ Locally,
$$\Psi(q^{i},y^{A},\dot{q}^{i},\dot{y}^{A},\mu_{i},\mu_{A},\gamma_{i},\gamma_{A})=(q^{i},y^{A},\gamma_{i},\gamma_{A}).$$
\begin{definition}

Define the \textbf{Legendre transform} associated with a
second-order nonholonomic system
$(\mathcal{D}^{(2)},\Sigma_{\mathcal{L}})$ as the map
$\mathbb{F}\mathcal{L}:\Sigma_{\mathcal{L}}\ra T^{*}\mathcal{D}$
given by $\mathbb{F}\mathcal{L}=\Psi\circ i_{\Sigma_{\mathcal{L}}}.$
In local coordinates, it is given by
$$\mathbb{F}\mathcal{L}(q^{i},y^{A},\dot{y}^{A},\gamma_{i})=\left(q^{i},y^{A},\gamma_{i},\frac{\partial\mathcal{L}}{\partial\dot{y}^{A}}\right).$$
\end{definition}

The following diagram summarizes the situation
$$
\xymatrix{
\Sigma_{\mathcal{L}} \ar@/_2pc/[rrrr]_{\mathbb{F}\mathcal{L}}  \ar@{^{(}->}[rr]^{i_{\Sigma_{\mathcal{L}}}} & \ & T^{*}T\mathcal{D} \ar[rr]^{\Psi} & \ & T^{*}\mathcal{D} \\
}
$$

\begin{definition}
We say that the second-order nonholonomic system
$(\mathcal{D}^{(2)},\Sigma_{\mathcal{L}})$ is \textbf{regular} if
$\mathbb{F}\mathcal{L}:\Sigma_{\mathcal{L}}\ra T^{*}\mathcal{D}$ is
a local diffeomorphism and \textbf{hyperregular} if
$\mathbb{F}\mathcal{L}$ is a global diffeomorphism.
\end{definition}

From the local expression of $\mathbb{F}\mathcal{L}$ we can observe
that from a direct application of the implicit function theorem we
have:
\begin{proposition}
The second-order nonholonomic system
$(\mathcal{D}^{(2)},\Sigma_{\mathcal{L}})$ determined by
$\mathcal{L}:\mathcal{D}^{(2)}\ra\R$ is regular if and only if the
matrix
$\displaystyle{\left(\frac{\partial^{2}\mathcal{L}}{\partial\dot{y}^{A}\partial\dot{y}^{B}}\right)}$
is non singular.
\end{proposition}

\begin{remark}
Observe that if the Lagrangian $\mathcal{L}:\mathcal{D}^{(2)}\ra\R$
is determined from an optimal control problem and its expression is
given by \eqref{lagrangiancontrol} then the regularity of the matrix
$\displaystyle{\left(\frac{\partial^{2}\mathcal{L}}{\partial\dot{y}^{A}\partial\dot{y}^{B}}\right)}$
 is equivalent to $$\det\left(\frac{\partial^{2}C}{\partial u^{A}\partial u^{B}}\right)\neq
 0$$ for the cost function.

\end{remark}

\subsection{Hamiltonian formalism}\


Assume that the system is regular. Then if we denote by
$p_{i}=\gamma_{i}$ and
$\displaystyle{p_{A}=\frac{\partial\mathcal{L}}{\partial\dot{y}^{A}}}$
we can write $\dot{y}^{A}=\dot{y}^{A}(q^{i},y^{A},p_{A}).$ Define
the \textit{Hamiltonian function}
$\mathcal{H}:T^{*}\mathcal{D}\ra\R$ by

$$\mathcal{H}(\alpha)=\langle\alpha,\pi_{T^{*}T\mathcal{D}}\mid_{\Sigma_{\mathcal{L}}}\left(\mathbb{F}\mathcal{L}^{-1}(\alpha)\right)\rangle-\mathcal{L}\left(\pi_{T^{*}T\mathcal{D}}\mid_{\Sigma_{\mathcal{L}}}\left(\mathbb{F}\mathcal{L}^{-1}(\alpha)\right)\right)$$
where $\alpha\in T^{*}\mathcal{D}$ is a one-form on $\mathcal{D},$
and
$\pi_{T^{*}T\mathcal{D}}\mid_{\Sigma_{\mathcal{L}}}:\Sigma_{\mathcal{L}}\ra\mathcal{D}^{(2)}$
is the projection locally given by
$\pi_{T^{*}T\mathcal{D}}\mid_{\Sigma_{\mathcal{L}}}(q^{i},y^{A},\dot{y}^{A},\gamma_{i})=(q^{i},y^{A},\dot{y}^{A}).$
Locally the Hamiltonian is given by
$$\mathcal{H}(q^{i},y^{A},p_{i},p_{A})=p_{A}\dot{y}^{A}(q^{i},y^{A},p_{A}))+p_{i}\rho_{A}^{i}y^{A}-\mathcal{L}(q^{i},y^{A},\dot{y}^{A}(q^{i},y^{A},p_{A})),$$
where we are using 
$$\mathbb{F}\mathcal{L}^{-1}(q^{i},y^{A},p_{i},p_{A})=\left(q^{i},y^{A},\rho_{A}^{i},\dot{y}^{A}(q^{i},y^{A},p_{A}),\frac{\partial\mathcal{L}}{\partial
q^{i}}-p_{j}\frac{\partial\rho_{A}^{j}}{\partial
q^{i}}y^{A},\frac{\partial\mathcal{L}}{\partial
y^{A}}-p_{j}\rho_{A}^{j},p_i,p_{A}\right).$$

Below we will see that the dynamics of the nonholonomic
optimal control problem is determined by the Hamiltonian system
given by the triple
$(T^{*}\mathcal{D},\omega_{\mathcal{D}},\mathcal{H})$ where
$\omega_{\mathcal{D}}$ is the standard symplectic $2-$form on
$T^{*}\mathcal{D}.$

The dynamics of the optimal control problem for the second-order
nonholonomic system is given by the symplectic hamiltonian dynamics
determined by the dynamical equation
\begin{equation}\label{dynamic}
i_{X_{\mathcal{H}}}\omega_{\mathcal{D}}=d\mathcal{H}.
\end{equation}  Therefore, if we consider the integral curves of $X_{\mathcal{H}},$ there are of the type $t\mapsto(\dot{q}^{i}(t),\dot{y}^{A}(t),\dot{p}_{i}(t),\dot{p}_{A}(t));$ the solutions of the nonholonomic Hamiltonian system is
specified by the Hamilton's equations on $T^{*}\mathcal{D}$

\begin{eqnarray*}
\dot{q}^{i}&=&\frac{\partial\mathcal{H}}{\partial p_{i}},\qquad\quad \dot{y}^{A}=\frac{\partial\mathcal{H}}{\partial p_{A}},\\
\dot{p}_{i}&=&-\frac{\partial\mathcal{H}}{\partial q^{i}},\qquad\dot{p}_{A}=-\frac{\partial\mathcal{H}}{\partial y^{A}};
\end{eqnarray*} that is,

\begin{eqnarray*}
\dot{q}^{i}&=&\rho_{A}^{i}y^{A},\\
\dot{p}_{i}&=&\frac{\partial\mathcal{L}}{\partial q^{i}}(q^{i},y^{A},\dot{y}^{A}(q^{i},y^{A},p_{A}))-p_j\frac{\partial\rho_{A}^{j}}{\partial q^{i}}y^{A},\\
\dot{p}_{A}&=&\frac{\partial\mathcal{L}}{\partial y^{A}}(q^{i},y^{A},\dot{y}^{A}(q^{i},y^{A},p_{A}))-p_{j}\rho_{A}^{j}.
\end{eqnarray*}

From equation \eqref{dynamic} it is clear that the flow preserves
 the symplectic $2-$form $\omega_{\mathcal{D}}.$ Moreover,
these equations are equivalent to equations given in
\eqref{eqELmulti} using the identification between the Lagrange
multipliers with the variables $p_i$ and the relation for
$\displaystyle{p_{A}}=\frac{\partial\mathcal{L}}{\partial\dot{y}^{A}}.$

\begin{remark}

We observe that in our formalism the optimal control
dynamics is deduced using a \textit{constrained variational
procedure} and equivalently it is possible to apply the 
Hamilton-Pontryagin's principle (see \cite{holm} for example), but,
in any case, this ``variational procedure'' implies the preservation
of the symplectic 2-form, and this is reflected in the Lagrangian
submanifold character. Moreover, in our case, under the regularity
condition, we have seen that the Lagrangian submanifold shows
that the system can be written as a Hamiltonian system (which is
obviously symplectic).

Additionally, we use the Lagrangian submanifold
$\Sigma_{\mathcal{L}}$ as a way to define intrinsically the
Hamiltonian side since we define the Legendre transformation using
the Lagrange submanifold $\Sigma_{\mathcal{L}}$. However there
exist other possibilities. For instance, in \cite{arnold2} (Section 4.2) the
authors  defined the corresponding momenta for a
vakonomic system. Using this procedure the momenta are locally expressed as follows

\begin{eqnarray*}
p_i&=&\frac{\partial \widetilde{\mathcal L}}{\partial \dot{q}^i}+\lambda^j \frac{\partial  f_j}{\partial \dot{q}^i}\\
p_A&=&\frac{\partial \widetilde{\mathcal L}}{\partial \dot{y}^A}+\lambda^j \frac{\partial  f^j}{\partial \dot{y}^A}
\end{eqnarray*} where $\widetilde{\mathcal L}$ is an arbitrary extension of ${\mathcal L}$ to $T{\mathcal D}$ and $f^j=\dot{q}^j-\rho^j_Ay^A=0$ are the constraint equations.  A simple computation shows that both are equivalent,
but  our
derivation is more intrinsic and geometric, that is, independent of
coordinates or extensions and without using Lagrange multipliers.



\end{remark}

\subsection{Example: continuously variable transmission (CVT)
(cont'd)}Now, we continue the example of the optimal control problem
for a continuously variable transmission  that we 
 considered in \hyperref[cvt]{Section $3.1$}. Recall that the
constraint distribution for the CVT is given by $\mathcal{D}\subset
T(\mathbb{S}^{1}\times\mathbb{S}^{1}\times\R)$

$$\mathcal{D}=\left\{\frac{1}{m}\frac{\partial}{\partial
x},(1-x)\frac{\partial}{\partial \theta_1}+x\frac{\partial}{\partial
\theta_2}\right\}.$$

The system is regular since
$$\det\left(\frac{\partial^{2}\mathcal{L}}{\partial\dot{y}^{A}\partial\dot{y}^{B}}\right)=\frac{(B(x))^2}{m^2}\neq
0$$ since $B(x)=J_1(1-x)^{2}+J_{2}x^{2}\neq 0$.

Denoting by
$(\theta_1,\theta_2,x,y_1,y_2,p_{\theta_1},p_{\theta_2},p_{x},p_1,p_2)$
local coordinates on $T^{*}\mathcal{D}$ the dynamic of the optimal
control problem for this nonholonomic system is determined by the
Hamiltonian function $\mathcal{H}:T^{*}\mathcal{D}\ra\R,$
\begin{eqnarray*}
\mathcal{H}(\theta,\theta_2,x,y_1,y_2,p_{\theta_1},p_{\theta_2},p_{x},p_1,p_2)&=&\frac{m^{2}p_{1}^{2}}{2}+\frac{p_2^{2}}{2(B(x))^{2}}+\frac{p_2A(x)y_1y_2}{mB(x)}+p_{\theta_1}(1-x)y_{2}\\
&+&p_{\theta_2}xy_{2}+p_{x}\frac{y_{1}}{m}.
\end{eqnarray*}
The corresponding Hamiltonian equations of motion are

\begin{eqnarray*}
\dot{y}_{1}&=&m^{2}p_{1},\quad\dot{p}_{\theta_1}=0,\\
\dot{y}_{2}&=&\frac{p_2}{(B(x))^{2}}+\frac{A(x)y_1y_2}{mB(x)},\quad\dot{p}_{\theta_2}=0,\\
\dot{p}_{x}&=&y_2(p_{\theta_1}-p_{\theta_2})-\frac{p_2y_1y_2((A(x))^{2}-J_1J_2)}{m(B(x))^{2}}-\frac{2p_{2}^{2}A(x)}{(B(x))^{3}},\\
\dot{p}_{1}&=&=-\frac{p_2A(x)y_2}{mB(x)}-\frac{p_x}{m},\quad\dot{p}_{2}=-\frac{p_2A(x)y_1}{mB(x)}-p_{\theta_1}(1-x)-p_{\theta_2}x.
\end{eqnarray*}

\subsection{Example: the Chaplygin sleigh (cont'd)}

In what follows, we continue the example of the optimal control
problem of the Chaplygin sleigh that we  began to study in
\hyperref[verticalcoin]{Section $3.2$}. Recall that the constraint
distribution is given
by $\mathcal{D}\subset TSE(2)$ where
$$\mathcal{D}=\left\{\frac{1}{J}\frac{\partial}{\partial\theta},\frac{\cos\theta}{m}\frac{\partial}{\partial x}+\frac{\sin\theta}{m}\frac{\partial}{\partial y}\right\}.$$

The system is regular since
$$\det\left(\frac{\partial^{2}\mathcal{L}}{\partial\dot{y}^{A}\partial\dot{y}^{B}}\right)=\frac{a^{4}}{J^{4}}\neq
0.$$

Denoting by $(x,y,\theta,y_1,y_2,p_x,p_y,p_{\theta},p_1,p_2)$ local
coordinates on $T^{*}\mathcal{D}$ the dynamics of the optimal control
problem for this nonholonomic system is determined by the
Hamiltonian function $\mathcal{H}:T^{*}\mathcal{D}\ra\R,$
$$
\mathcal{H}(x,y,\theta,y_1,y_2,p_x,p_y,p_{\theta},p_1,p_2)=\frac{J^{2}}{2b^2}p_{1}^{2}+\frac{m^2}{2}p_2^{2}+p_{x}\frac{\cos\theta}{m}y_{2}+\frac{p_{\theta}}{J}y_{1}+p_{y}\frac{\sin\theta}{m}y_{2}.$$

The Hamiltonian equations of motion are

\begin{eqnarray*}
\dot{y}_{1}&=&\frac{J^{2}p_{1}}{b^2},\quad\dot{y}_{2}=m^{2}p_{2},\quad\dot{p}_{x}=0,\quad \dot{p}_{y}=0,\\
\dot{p}_{\theta}&=&p_{x}\frac{\sin\theta}{m}y_{2}-p_{y}\frac{\cos\theta}{m}y_{2},\\
\dot{p}_{1}&=&=-\frac{p_{\theta}}{J},\quad\dot{p}_{2}=-p_{x}\frac{\cos\theta}{m}-p_{y}\frac{\sin\theta}{m}.
\end{eqnarray*}

 Integrating the equations $\dot{p}_{x}=0$ and $\dot{p}_{y}=0$
as $p_x=c_1$ and $p_y=c_2$ where $c_1$ and $c_2$ are constants the
 system of differential equations becomes

\begin{eqnarray*}
\dot{y}_{1}&=&\frac{J^{2}p_{1}}{b^2},\quad\dot{p}_{\theta}=c_{1}\frac{\sin\theta}{m}y_{2}-c_{2}\frac{\cos\theta}{m}y_{2},\\
\dot{y}_{2}&=&m^{2}p_{2},\quad\dot{p}_{1}=-\frac{p_{\theta}}{J},\quad\dot{p}_{2}=-c_{1}\frac{\cos\theta}{m}-c_{2}\frac{\sin\theta}{m}.
\end{eqnarray*}

Differentiating  $\dot{y}_{1}$ and
$\dot{y}_{2}$ and substituting we obtain
$$
\frac{\dddot{y}_{1}}{J}=\frac{y_{2}}{mb^{2}}\left(c_{2}\cos\theta-c_{1}\sin\theta\right),\quad
\ddot{y}_{2}=-m(c_{1}\cos\theta+c_2\sin\theta),
$$  as in the Lagrangian setting.

Observe that in the case of motion planing for obstacle avoidance
the Hamiltonian function $\mathcal{H}:T^{*}\mathcal{D}\to\R$ is
given by
\begin{eqnarray*}\mathcal{H}(x,y,\theta,y_1,y_2,p_x,p_y,p_{\theta},p_1,p_2)&=&\frac{J^{2}}{2b^2}p_{1}^{2}+\frac{m^2}{2}p_2^{2}+p_{x}\frac{\cos\theta}{m}y_{2}+\frac{p_{\theta}}{J}y_{1}+p_{y}\frac{\sin\theta}{m}y_{2}\\
&-&\frac{\kappa}{2(x-x_C)^{2}+2(y-y_C)^{2}},
\end{eqnarray*}

and the resulting dynamical equations are
\begin{eqnarray*}
\dot{y}_{1}&=&\frac{J^{2}p_{1}}{b^2},\quad\dot{y}_{2}=m^{2}p_{2}\;,\quad\dot{p}_{x}=\frac{\kappa(x-x_C)}{((x-x_C)^{2}+(y-y_C)^{2})^2},\\
\dot{p}_{y}&=&\frac{\kappa(y-y_C)}{((x-x_C)^{2}+(y-y_C)^{2})^2}\; ,\quad\dot{p}_{\theta}=p_{x}\frac{\sin\theta}{m}y_{2}-p_{y}\frac{\cos\theta}{m}y_{2},\\
\dot{p}_{1}&=&=-\frac{p_{\theta}}{J},\quad\dot{p}_{2}=-p_{x}\frac{\cos\theta}{m}-p_{y}\frac{\sin\theta}{m}.
\end{eqnarray*}

\section{Conclusions and future research}\label{conclusiones}
In this section we summarize the
contributions of our work and discuss future research.
\subsection{Conclusions:}
In this paper we study optimal control problems for a class of
nonholonomic mechanical systems. We have given  a geometrical derivation of the equations of motion of a nonholonomic optimal control problem as a constrained problem on the tangent space to the constraint distribution. We have seen how the dynamics of the optimal control problem can be completely described by a Lagrangian submanifold of an appropriate cotangent bundle and under some mild regularity conditions we have derived the the equations of motion for the nonholonomic optimal control problem as a classical set of Hamilton's equations on the cotangent bundle of the constraint distribution.  We have introduced the notion of Legendre transformation in this context to establish the relationship between  the Lagrangian and Hamiltonian dynamics. We applied our techniques to different examples: optimal control of a continuously variable transmission, Chaplygin sleigh and to optimal planning for obstacle avoidance problems.

\subsection{Future research: Construction of geometric and variational integrators for optimal control problems of nonholonomic mechanical systems.}
In this paper we have seen that an  optimal control problem of a
nonholonomic system may be viewed  as a Hamiltonian system on
$T^{*}\mathcal{D}$. One can thus use standard methods for
symplectic integration such as symplectic Runge-Kutta methods,
collocation methods, St$\ddot{\hbox{o}}$rmer-Verlet, symplectic
Euler methods, etc.; developed and studied in \cite{EqAlge},
\cite{LeRe1}, \cite{LeRe2}, \cite{SS}, \cite{SSerna}, e.g., to
simulate nonholonomic optimal control problems. 

Also, we would like to build variational integrators as an
alternative way to construct integration schemes for these kinds of
optimal control problems following the results given in
\hyperref[section3]{Section 3}. Recall that in the
continuous case we have considered a Lagrangian function
$\mathcal{L}:\mathcal{D}^{(2)}\to\R$. Since the space
$\mathcal{D}^{(2)}$ is a subset of $T\mathcal{D}$ we can discretize
the tangent bundle $T\mathcal{D}$ by the cartesian product
$\mathcal{D}\times\mathcal{D}$. Therefore, our discrete variational
approach for optimal control problems of nonholonomic mechanical
systems will be determined by the construction of a discrete
Lagrangian $\mathcal{L}_{d}:\mathcal{D}_{d}^{(2)}\to\R$ where
$\mathcal{D}^{(2)}_{d}$ is the subset of
$\mathcal{D}\times\mathcal{D}$ locally determined by imposing the
discretization of the constraint $\dot{q}^{i}=\rho_{A}^{i}(q)y^{A}$,
for instance we can consider
$$\mathcal{D}_{d}^{(2)}=\left\{(q_0^{i},y_{0}^{A},q_{1}^{i},y^{A}_{1})\in
\mathcal{D}\times\mathcal{D}\bigg{|}\frac{q_{1}^{i}-q_0^{i}}{h}=\rho_{A}^{i}\left(\frac{q_0^{i}+q_0^{i}}{2}\right)\left(\frac{y_{0}^{A}+y_{1}^{A}}{2}\right)\right\}.$$

Now the system is in a form appropriate for the application of discrete
variational methods for constrained systems (see \cite{MMS} and references therein).
\subsection*{Acknowledgment} We wish to thanks Klas Modin and Olivier Verdier the permission to use their graphical illustration and description of the Continuously Variable Transmission Gearbox.

\end{document}